\documentclass[reqno]{amsart}
\usepackage{hyperref}
\usepackage{dsfont} 
\usepackage{amssymb} 
\usepackage{extarrows}
\usepackage{mathrsfs} 

\AtBeginDocument{{\noindent\small
Vol. 2018 (2018), No. , pp. .\newline
ISSN: . URL:}
\vspace{8mm}}

\begin{document}
\title[\hfilneg EJDE-2018/08\hfil Invariant foliations for stochastic dynamical systems]
{Invariant foliations for stochastic dynamical systems with  multiplicative stable L\'{e}vy noise}

\author[Ying Chao, Pingyuan Wei, Shenglan Yuan \hfil EJDE-2018/08 \hfilneg]
{Ying Chao, Pingyuan Wei$^\ast$, Shenglan Yuan}

\address{Ying Chao \newline
School of Mathematics and Statistics \& Center for Mathematical Sciences, Hubei Key Laboratory of Engineering Modeling and Scientific Computing,
Huazhong University of Sciences and Technology,
Wuhan 430074, P.R. China}
\email{yingchao1993@hust.edu.cn}

\address{$^\ast$Corresponding Author: Pingyuan Wei \newline
School of Mathematics and Statistics \& Center for Mathematical Sciences, Hubei Key Laboratory of Engineering Modeling and Scientific Computing,
Huazhong University of Sciences and Technology,
Wuhan 430074, P.R. China}
\email{weipingyuan@hust.edu.cn}

\address{Shenglan Yuan \newline
School of Mathematics and Statistics \& Center for Mathematical Sciences, Hubei Key Laboratory of Engineering Modeling and Scientific Computing,
Huazhong University of Sciences and Technology,
Wuhan 430074, P.R. China}
\email{shenglanyuan@hust.edu.cn}

\dedicatory{Communicated by }

\thanks{Submitted . Published .}
\keywords{Stochastic differential equations; random dynamical systems;
\hfill\break\indent invariant foliations; invariant manifolds;  geometric structure}

\begin{abstract}
This work deals with the dynamics of a class of stochastic dynamical systems with a multiplicative non-Gaussian L\'{e}vy noise. We first establish the existence of stable and unstable foliations for this kind of system via the Lyapunov-Perron method. Then we examine the geometric structure of the invariant foliations, and their relation with invariant manifolds.  Finally, we illustrate our results in an example.
\end{abstract}

\maketitle
\newtheorem{theorem}{Theorem}[section]
\newtheorem{example}[theorem]{Example}
\newtheorem{Theorem}{Theorem}[section]
\newtheorem{Definition}[Theorem]{Definition}
\newtheorem{Proposition}[Theorem]{Proposition}
\newtheorem{Lemma}{Lemma}[section]
\newtheorem{Example}[Theorem]{Example}
\newtheorem{Claim}[Theorem]{Claim}
\newtheorem{Corollary}{Corollary}[section]
\newtheorem{Remark}{Remark}[section]

\section{Introduction} \label{Intr}

\noindent Invariant foliations, as well as invariant manifolds, are geometric structures in state space for describing and understanding the dynamics of nonlinear dynamical systems (\cite{BW, CDLS,DLS-2003, LS, MZZ}).
An invariant foliation is about describing sets (called fibers) in state space with certain dynamical properties. A fiber consists of all those points starting from which the dynamical orbits are exponentially approaching each other, in forward time (stable foliation) or backward time (unstable foliation). Both stable and unstable fibers are building blocks for   dynamical systems, as they carry specific dynamical information.
The stable and unstable foliations for deterministic systems have been investigated by various authors \cite{BLZ,BLZ1,CHT,CLL,Henry}.

\par During the last two decades, there has been various studies on invariant foliations and invariant manifolds for stochastic differential equations (SDEs). Lu and Schmalfuss \cite{LS08} proved the existence of random invariant foliations for infinite dimensional stochastic dynamical systems. Moreover, Sun et al. \cite{SKD} provided an approximation method of invariant foliations for dynamical systems with small noisy perturbations via asymptotic analysis. Subsequently, Chen et al. \cite{Cdz} further studied the slow foliation of a multiscale (slow-fast) stochastic evolutionary system, eliminating the fast variables for this system. Most of these works were for stochastic systems with Gaussian noise, i.e., Brownian noise.

However, in applications of biological and physical fields, noise appeared in the complex systems are often non-Gaussian rather than Gaussian \cite{Zheng2016, WangHui, BSW,Wo,IM}.
Note that the slow manifolds of a class of slow-fast stochastic dynamical systems with non-Gaussian additive type noise and its approximation have been considered by Yuan et al. \cite{YD}. Kummel \cite{Kummel} studied invariant manifolds of finite dimensional stochastic systems with multiplicative noise.
It is now the time to consider invariant foliations for stochastic dynamical systems with non-Gaussian noise.

 In this paper, we are concerned with invariant foliations for stochastic systems in case of non-Gaussian L\'evy noise and their relationship with invariant manifolds.
\par

Consider the following nonlinear stochastic dynamical system with linear multiplicative $\alpha$-stable L\'{e}vy noise
\begin{eqnarray}
\frac{dx}{dt}&=&Ax+f(x,y)+x\diamond\dot{L_t^\alpha},~~~~~~\quad \hbox{in}\quad \mathbb{R}^n, \label{Equation-s}\\
\frac{dy}{dt}&=&By+g(x,y)+y\diamond\dot{L_t^\alpha},~~~~~~\quad \hbox{in}\quad \mathbb{R}^m, \label{Equation-f}
\end{eqnarray}
where $\diamond$ denotes Marcus differential \cite{Ap,Kummel}. The operators A and B are generators of $C_0$-semigroups satisfying an exponential dichotomy condition. Nonlinearities $f$ and $g$ are Lipschitz continuous functions with $f(0,0)=0$, $g(0,0)=0$. The stochastic process ${L_t^\alpha}$ is a scalar, two-sided symmetric $\alpha$-stable L\'{e}vy process with index of the stability ${1<\alpha<2}$ \cite{Ap,Duan}. The precise conditions
on these quantities will be specified in Section $3$.
\par
It is worthy mentioning that as Marcus SDEs preserve certain physical quantities such as energy, they are often appropriate models in engineering and physical applications \cite{Sdf}.  The linear multiplicative noise appears in the cases where noise fluctuates in proportion to the system state, as in some geophysical systems and fulid systems. The wellposedness of mild solutions for this kind of stochastic differential equations with non-Gaussian L\'evy noise is known \cite{FK,Kunita, Ap, PZ}.

\par
 In order to provide a geometric visualization for the state space of dynamical system \eqref{Equation-s}-\eqref{Equation-f} via invariant foliations in the similar sprit as in \cite{FLD,Kummel}, and to explore its geometry structure, we first introduce a random transformation based on the L\'{e}vy-type Ornstein-Uhlenbeck process to convert a Marcus SDE into a conjugated random differential equation (RDE) which easily generates a random dynamical system. Then we prove that, under appropriate conditions, an unstable foliation can be constructed as a graph of a Lipschitz continuous map via the Lyapunov-Perron method \cite{Boxler,FLD}. After that, by the inverse transformation, we can obtain the unstable foliation for the original stochastic system. Furthermore, we shall analyze the geometric structure of the unstable foliation  and verify that the unstable manifold is one fiber of the unstable foliation. There are similar conclusions about the stable foliation.

\par
This paper is arranged as follows. In Section 2, we offer a brief summary of basic concepts in random dynamical systems and present a special but very important metric dynamical system represented by a L\'{e}vy process with two-sided time. Subsequently, Marcus canonical differential equations with L\'{e}vy motions are discussed. Our framework is presented in Section 3. In Section 4, we show the existence of unstable foliation (Theorem \ref{Th4.1}), examine its geometric structure and illustrate a link with unstable manifold (Theorem \ref{Th4.3}). The same results on the stable foliation for \eqref{Equation-s}-\eqref{Equation-f} are given in Theorem \ref{Th4.4}. Finally, Section 5 is devoted to an illustrative example.
\par

\renewcommand{\theequation}{\thesection.\arabic{equation}}
\setcounter{equation}{0}

\section{Preliminaries }
\quad\; We now recall some preliminary concepts in random dynamical systems \cite{Ap,Arnold,Kummel}. Then we discuss differential equations driven by L\'{e}vy noise.
\par

\subsection{Random dynamical systems }

\quad \;
Let us recall an appropriate model for noise. \par
\begin{Definition}(Metric dynamical system)\label{MDS}
 Given a probability space $(\Omega,\mathcal{F}, \mathds{P})$ and a flow  $\theta$=$\{\theta_t\}_{t\in \mathbb{R}}$  on $\Omega$ defined as a mapping
 $$
 \theta: \mathbb{R}\times \Omega\mapsto \Omega,
 $$
 which satisfies
\par
$\bullet$ $\theta_0=id  (indentity)~~on~~ \Omega; $
\par
$\bullet$ $\theta_{t_1}\theta_{t_2}=\theta_{t_1+t_2}~~~~for~~all~~ t_1, t_2 \in \mathbb{R}; $
\par
$\bullet$ the mapping $(t,\omega)\mapsto \theta_t\omega$ is $(\mathcal{B}(\mathbb{R})\otimes\mathcal{F},
\mathcal{F})-measurable$, where $\mathcal{B}(\mathbb{R})$ is the collection of Borels sets on the real line $\mathbb{R}$.
\par
\noindent In addition, the probability measure $\mathds{P}$ is
assumed to be ergodic with respect to $\{\theta_t\}_{t\in
\mathbb{R}}$. Then the quadruple $\Theta: =(\Omega, \mathcal{F}, \mathds{P}, \theta)$
is called a metric dynamical system.
\end{Definition}

\par
For our applications, we will consider a canonical sample space for two-sided L\'{e}vy process. Let $\Omega={\mathcal D}(\mathbb{R},\mathbb{R}^d)$ be
 the space of c\`{a}dl\`{a}g functions (i.e., continuous on the right and have limits on the left) taking zero value at $t = 0$ defined on $\mathbb{R}$ and taken values in $\mathbb{R}^d$. $\mathcal D(\mathbb{R}, \mathbb{R}^d)$ is not separable if we use the usual compact-open metric. To make space $\mathcal D(\mathbb{R},\mathbb{R}^d)$ complete and separable, a Skorokhod's topology generated by the Skorokhod's metric $d_\mathbb{R}$ is equipped \cite{Billingsley,Situ}. For functions $\omega_1,\omega_2\in \mathcal D(\mathbb{R}, \mathbb{R}^d), d_\mathbb{R}(\omega_1, \omega_2)$ is defined as
 $$d_\mathbb{R}(\omega_1,\omega_2):=\sum_{n=1}^\infty\frac{1}{2^n}(1\wedge d_n(\omega_1^n,\omega_2^n)), $$
 where $\omega_1^n(t):=f_n(t)\omega_1(t),\omega_2^n(t):=f_n(t)\omega_2(t) $ with
 $$f_n(t)=\left\{
\begin{aligned}
&1\quad ~~~~~if ~~~~|t|\leq n-1;\\
&n-t  \quad if ~~~n-1\leq|t|\leq n; \\
&0 \quad ~~~~~if ~~~ |t|\geq n.
\end{aligned}
\right.
$$
and
$$d_n(\omega_1^n,\omega_2^n):=\inf_{\lambda\in\Lambda}\{\sup_{-n\leq s<t\leq n}|\ln\frac{\lambda(t)-\lambda(s)}{t-s}\mid\vee \sup_{-n\leq t\leq n}|\omega_1(t)-\omega_2(\lambda(t))|\},$$
where $\Lambda:=\{\lambda:\mathbb{R}\to \mathbb{R};\lambda   \hbox{ is strictly increasing},
\lim_{t \to -\infty}\lambda(t)=-\infty,\lim_{t \to +\infty}\lambda(t)=+\infty\}.$
We denote by $\mathcal{F}$:= $\mathcal{B}(\mathcal D(\mathbb{R},\mathbb{R}^d))$ the associated Borel $\sigma$-algebra. On this set, measurable flow $\theta = \{\theta_t\}_{t\in \mathbb{R}}$ is defined
by the shifts
$$\theta_t\omega =\omega(\cdot+t)-\omega(t),\quad \omega\in\Omega,\quad t\in \mathbb{R}.$$
Let $\mathds{P}$ be the probability measure on $\mathcal{F},$ which is given by
the distribution of a two-sided L\'{e}vy motion with path in $\mathcal D(\mathbb{R},\mathbb{R}^d)$. Note that $\mathds{P}$ is ergodic with respect to $\theta_t$; see \cite[Appendix A]{Arnold}.
Thus $(\mathcal D(\mathbb{R},\mathbb{R}^d),\mathcal{B}(\mathcal D(\mathbb{R},\mathbb{R}^d)),\mathds{P},\{\theta_t\}_{t\in \mathbb{R}})$ is a metric dynamical syatem. Later on we will consider, instead of the whole $\mathcal D(\mathbb{R},\mathbb{R}^d)$, a $\{\theta_t\}_{t\in \mathbb{R}}$-invariant subset $\Omega\subset \mathcal D(\mathbb{R},\mathbb{R}^d)$ of $\mathds{P}$-measure one as well as the trace $\sigma$-algebra $\mathcal{F}$ of
$\mathcal{B}(\mathcal D(\mathbb{R},\mathbb{R}^d))$ with respect to $\Omega$. Review that a set $\Omega$ is called {$\{\theta_t\}_{t\in \mathbb{R}}$-invariant if $\theta_t\Omega = \Omega$ for $t \in \mathbb{R}$ \cite[ Page545] {Arnold}. On $\mathcal{F}$, we will consider the restriction of the measure $\mathds{P}$
and still denote it by $\mathds{P}$. In our set, we consider scalar L\'{e}vy motion, i.e., d=1.
\begin{Definition}(Random dynamical system (RDS))\label{RDS}
 A measurable random dynamical system on a measurable space $(H,\mathcal{B}(H))$ over the metric dynamical system $(\Omega,\mathcal{F}, \mathds{P}, \theta)$ is given by a mapping:
$$
\varphi:    \mathbb{R}\times\Omega\times H \mapsto H,
$$
with the following properties:
\par $\bullet$ $\varphi$ is jointly $(\mathcal{B}(\mathbb{R})\otimes \mathcal{F} \otimes
\mathcal{B}(H), \mathcal{B}(H))$-measurable;
\par
$\bullet$ the mapping $\varphi(t,\omega)$:=$\varphi(t,\omega,\cdot)$: $H\mapsto H$ form a cocycle over $\theta(\cdot)$, that is:
$$
\begin{array}{l}
\varphi(0, \omega, x) = x,\\
\varphi(t_1 + t_2, \omega, x) = \varphi(t_2, \theta_{t_1}\omega, \varphi(t_1,
\omega, x)),
\end{array}
$$
for each $t_1, t_2 \in \mathbb{R}$, $\omega \in \Omega$ and $x\in H$. In this paper, we take H=$\mathbb{R}^{n+m}=\mathbb{R}^{n}\times\mathbb{R}^{m}$.
\end{Definition}
\par


Generally speaking, a stable foliation or an unstable foliation is composed of stable fibers or unstable fibers which are certain sets in the state space carrying specific  dynamical information. More precisely, a stable fiber or an unstable fiber of a foliation is defined as follows \cite{Cdz,CLL}:
\begin{Definition}(Stable and unstable fiber)\label{SUF}
\par
\noindent(i) $\mathcal{W}_{\gamma s}(x, \omega)$ is called a $\gamma$-stable fiber passing through $x\in H$ with $\gamma\in
\mathbb{R}^-$,  if $\|\varphi(t,\omega,x)-\varphi(t, \omega, \hat{x})\|_H=O(e^{\gamma
t}), \forall\; \omega\in \Omega$ as $t\to +\infty$ for all  $x,
\hat{x}\in \mathcal{W}_{\gamma s}$.
\par
\noindent(ii) $\mathcal{W}_{\eta u}(x, \omega)$ is called a
$\eta$-unstable fiber passing through $x\in H$ with $\eta\in
\mathbb{R}^+$, if $\|\varphi(t, \omega,x)-\varphi(t, \omega,\hat{x})\|_H=O(e^{\eta t}),
\forall\; \omega\in \Omega$ as $t\to -\infty$ for all $x, \hat{x}\in \mathcal{W}_{\eta u}$.
\end{Definition}
\par

 From the proceeding  definition, we see that a stable fiber or an unstable fiber is the set of all those points passing through which the dynamical trajectories can approach each other exponentially, in forward time or backward time, respectively. In fact, we can replace $O(e^{\gamma t})$ by $O(e^{p\gamma t})$  with $0<p\leq1$ as we will show,  without affecting the property of exponential approximation. In addition, we say a foliation is \emph{invariant }if the random
dynamical system $\varphi$ maps one fiber to another fiber in the
following sense
$$
\varphi(t, \omega, \mathcal{W}_\eta(x, \omega))\subset
\mathcal{W}_\eta(\varphi(t,\omega, x), \theta_t\omega).
$$

\subsection{Marcus canonical stochastic differential equations with  L\'{e}vy motions}

\quad\;
In the present, we consider a special but very useful class of scalar L\'{e}vy motions, i.e., the symmetric $\alpha$-stable L\'{e}vy motions $(1<\alpha<2)$ with drift zero, diffusion $d>0$ and L\'{e}vy measure $\nu_\alpha(du)=c_\alpha\frac{du}{|u|^{1+\alpha}}$ where $c_\alpha=\frac{\alpha}{2^{1-\alpha}\sqrt\pi}\frac{\Gamma(\frac{1+\alpha}{2})}{\Gamma(1-\frac{\alpha}{2})}$. Here  $\Gamma$ is Gamma function. For more definition, see \cite{Ap,Sato}.
\par
Initially Marcus canonical differential equations with point process as the driving process were discussed by Marcus in \cite{Marcus}. Subsequently, Kurtz $et$ $al.$ \cite{KPP} generalized the driving process. For a scalar symmetric L\'{e}vy motion $L_t^\alpha$ mentioned above, the precise definition is given by
$$dx(t)=b(x(t))dt+\sigma(x(t-))\diamond dL_t^\alpha$$
where $\diamond$ denotes the Marcus integral, i.e.,
\begin{eqnarray}
dx(t)=b(x(t))dt+\sigma(x(t-))\circ dL^{\alpha,c}(t)+\sigma(x(t-))d L^{\alpha,d}(t)\notag\\
+\sum_{0<s\leq t}[\psi(x(s-),\Delta L_s^\alpha)-x(s-)-\sigma(x(s-))\Delta L_s^\alpha],\notag
\end{eqnarray}
where $L^{\alpha,c}(t),L^{\alpha,d}(t)$ are the continuous and discontinuous parts of $L_t^\alpha$ respectively, $\circ$ denotes the Stratonovich integral. Moreover, $\psi(x(s-),\Delta L_s^\alpha)=\varsigma(\Delta L_t^\alpha\sigma;x(t-),1)$ satisfies equation:
$$\frac{d\varsigma(\sigma;v,t)}{dt}=\sigma[\varsigma(\sigma;v,t)],\quad \varsigma(\sigma;v,0)=v.$$
Appropriate conditions for coefficients $b$ and $\sigma$ as we will give later can ensure the existence and uniqueness of solution of the Marcus canonical equation, and then it defines a stochastic flow or cocycle so that RDS methods can be applied. For more details, see \cite{FK,Kunita}.
\par
Here, the reason for taking $1<\alpha<2$ is to ensure that Lemma \ref{L3.1} holds. In fact, the index of stability can take values in $(0,2)$. When $\alpha=2$, it reduces to the well-known Brownian motion.

\section{Framework}
\label{slow888}

\quad\;
For the system \eqref{Equation-s}-\eqref{Equation-f}, let $|\cdot|$ denotes the Euclidean norm. To construct the unstable foliation of system, we need to introduce the following hypotheses.
\par

{\bf A1 }(Exponential dichotomy condition): The linear operator $A$ be the generator of a $C_0$-semigroup $e^{At}$ on $\mathbb{R}^n$ satisfying
$$
|e^{At}x|\leq e^{a t}|x|,\quad
\hbox{for}\quad t\leq 0.
$$
Moreover, the linear operator $B$ is the generator of a $C_0$-semigroup $e^{Bt}$ on $\mathbb{R}^m$
satisfying
$$
|e^{Bt}y|\leq e^{b t}|y|,\quad
\hbox{for}\quad t\geq 0,
$$
where $b<0<a$.
\par

{\bf A2} (Lipschitz condition): The interactions functions
$$
\begin{array}{l}
f: \mathbb{R}^n\times \mathbb{R}^m \longrightarrow \mathbb{R}^n,\\
g: \mathbb{R}^n\times \mathbb{R}^m \longrightarrow \mathbb{R}^m,
\end{array}
$$
  are Lipschitz continuous with $f(0,0)=0$ and $g(0,0)=0$, i.e., there exists a positive
constant $K$ such that for all $(x_i,y_i) \in \mathbb{R}^n\times \mathbb{R}^m$, $i=1, 2,$
$$
|f(x_1, y_1)-f(x_2, y_2)|\leq
K(|x_1-x_2|+|y_1-y_2|),
$$
$$
|g(x_1, y_1)-g(x_2, y_2)|\leq
K(|x_1-x_2|+|y_1-y_2|).
$$

\par

\noindent Note that if $f$ and $g$ are locally Lipschitz, following the analysis in this paper, we also get invariant foliation in a neighborhood of (0, 0).\par
As in references \cite{DLS-2003,DLS-2004}, we are going to verify that stochastic system \eqref{Equation-s}-\eqref{Equation-f} can be transformed into the random differential system which is described by differential equations with random coefficients.
For this purpose, we consider a Langevin equation:
\begin{equation}\label{OU}
dz=-z dt + dL_t^\alpha.
\end{equation}
A solution of this equation is usually called a L\'{e}vy-type Ornstein-Uhlenbeck process. The properties of its stationary solution can be characterized by the following lemma in the same sprit of the case of Brownian noise, refer to \cite{DLS-2003,LDLK}.
\par


\begin{Lemma}\label{L3.1}
Let $L_t^\alpha$ be a two-sided  scalar symmetric $\alpha$-stable L\'{e}vy motion   with ${1<\alpha<2}$. Then
\par
\noindent(i) there exists a $\{\theta_t\}_{t\in \mathbb{R}}$-invariant set $\Omega\subset D(\mathbb{R},\mathbb{R}^d))$ of full measure with sublinear growth:
$$\lim_{t \to \pm\infty}\frac{\omega(t)}{t}=0, \quad\omega \in \Omega$$ of $\mathds{P}$-measure one.
\par
\noindent(ii) for $\omega \in \Omega$, the random variable
$$ z(\omega)=-\int_{-\infty}^0e^{\tau}\omega(\tau)d\tau$$ exists and generates a unique c\`{a}dl\`{a}g stationary solution of (\eqref{OU}) given by
$$ z(\theta_t\omega)=-\int_{-\infty}^0e^{\tau}\theta_t\omega(\tau)d\tau= -\int_{-\infty}^0e^{\tau}\omega(\tau+t)d\tau+\omega(t).  $$
\par
\noindent(iii) in particular,
$$\lim_{t \to \pm\infty}\frac{\left|z(\theta_t\omega)\right|}{\left| t \right|}=0, \quad\omega \in \Omega.$$
\par
\noindent(iv) in addition,
$$\lim_{t \to \pm\infty}\frac{1}{t}\int_{0}^tz(\theta_\tau\omega)d\tau=0, \quad\omega \in \Omega.$$
\end{Lemma}
\par
\noindent{\bf Proof.} (i) The $\int_{|x|>1}|x|\nu_\alpha(dx)$ is finite due to $\alpha$ has a value between 1 and 2 (\cite[Page 80]{Sato}). In addition, $E|L_1^\alpha|<\infty$ and $EL_1^\alpha=0$, by the properties of moments for L\'{e}vy process \cite[Page 163]{Sato}. Thus the assertion is obtained from the strong law of large numbers for L\'{e}vy process (see \cite[Page 246]{Sato}).
\par
\noindent (ii) The existence of the integral on the right hand side for $\omega\in \Omega_2$ follows from the fact that the sample paths of an $\alpha$-stable L\'{e}vy motion satisfy $\limsup_{t\to\infty}t^{-\frac{1}{\eta}}L_t^{\alpha,\ast}=0\quad a.s.\quad or=\infty\quad a.s.$ according to whether $\eta<\alpha$ or $\eta>\alpha$, respectively, where $L_t^{\alpha,\ast}=\sup_{0\leq s\leq t}|L_s^\alpha|$. For the remaining part we refer to \cite[Page 216 and 311]{Ap}.
\par
\noindent (iii) Based on the above facts, for $\frac{1}{\alpha}<\delta<1$ and $\omega\in \Omega_2$, there exists a constant $C_{\delta,\omega}>0$ such that $|\omega(\tau+t)|\leq C_{\delta,\omega}+|\tau|^\delta+|t|^\delta.$ Thus, $\lim_{t\to\pm\infty}|-\int_{-\infty}^0e^{\tau}\omega(\tau+t)d\tau|=0$, (iii) is proven.
\par
\noindent (iv) Since $L_t^\alpha$ is symmetric $\alpha$-stable, we can prove that $z(\omega)$ is also symmetric $\alpha$-stable, and $Ez(\omega)=0$. Thus, by the ergodic theorem we obtain (iv) for $\omega\in \Omega_3\in\mathcal{B}(D(\mathbb{R},\mathbb{R}^d))$. This set $\Omega_3$ is also {$\{\theta_t\}_{t\in \mathbb{R}}$-invariant. Then we set $\Omega:=\Omega_1\cap\Omega_2\cap\Omega_3$. The proof is complete. \hfill{$\blacksquare$}
\par
From now on, we replace $\mathcal{B}(D(\mathbb{R},\mathbb{R}^d))$ by
$$ \mathcal{F}=\{\Omega\cap A,A\in\mathcal{B}(D(\mathbb{R},\mathbb{R}^d))\}$$
for $\Omega$ given in Lemma \ref{L3.1}. Probability measure is the restriction of the original measure to this new $\sigma$-algebra, we still denote it by $\mathds{P}$.
\par
Define the random transformation
\begin{equation}\label{Transformation}
\binom{\hat{x}}{\hat{y}}=T(\omega,x,y):=\binom{xe^{-z(\omega)}}{ye^{-z(\omega)}}
\end{equation}
According to \cite{DLS-2004,ERV}, marcus canonical integral satisfies the usual chain rules, thus, the $(\hat{x}(t), \hat{y}(t))=T(\theta_t\omega, x(t), y(t))$ satisfies
the following conjugated random differential equations:
\begin{eqnarray}
\frac{d\hat{x}}{dt}&=&A\hat{x}+F(\hat{x}, \hat{y}, \theta_t\omega)+z(\theta_t\omega)\hat{x},  \label{RE-s} \\
\frac{d\hat{y}}{dt}&=&B\hat{y}+G(\hat{x},\hat{y},\theta_t\omega)+z(\theta_t\omega)\hat{y}, \label{RE-f}
\end{eqnarray}
where $$F(\hat{x},\hat{y},\theta_t\omega):=e^{-z(\theta_t\omega)}f(e^{z(\theta_t\omega)}\hat{x},e^{z(\theta_t\omega)}\hat{y}), $$
$$G(\hat{x},\hat{y},\theta_t\omega):=e^{-z(\theta_t\omega)}g(e^{z(\theta_t\omega)}\hat{x},e^{z(\theta_t\omega)}\hat{y}) $$
and $z(\theta_t\omega)$ is the c\`{a}dl\`{a}g stationary solution of \eqref{OU} given in lemma 3.1. We can see that functions $F$ and $G$ also satisfy the Lipschitz condition with the same Lipschitz constant $K$. Here, it is worth noting that although $x,y$ are only c\`{a}dl\`{a}g in time, the solution $\hat{x}(t),\hat{y}(t)$ are the product of two c\`{a}dl\`{a}g functions, and are actually continuous in time. And $\dot{\hat{x}}$ denote $\frac{d}{dt^{+}}\hat{x}:=\lim_{h\downarrow0^+}\frac{\hat{x}(t+h)-\hat{x}(t)}{h}$, i.e., the right derivations of $\hat{x}$ with respect to $t$. Therein, the state space for this new system is still $\mathbb{R}^{n+m} =\mathbb{R}^{n} \times \mathbb{R}^{m}$.
\par
Let $\hat{z}(t,\omega,\hat{z}_0):=(\hat{x}(t,\omega,(\hat{x}_0,\hat{y}_0)),\hat{y}(t,\omega,(\hat{x}_0,\hat{y}_0)))$ be the mild solution of the \eqref{RE-s}-\eqref{RE-f} with initial values $(\hat{x}(0),\hat{y}(0))=(\hat{x}_0,\hat{y}_0):=\hat{z}_0$ in the sense of Carath\'{e}odory \cite{ERV}. Then, the solution operator of \eqref{RE-s}-\eqref{RE-f},
$$\varphi(t,\omega,(\hat{x}_0,\hat{y}_0))=(\hat{x}(t,\omega,(\hat{x}_0,\hat{y}_0)),\hat{y}(t,\omega,(\hat{x}_0,\hat{y}_0)))$$
generate a random dynamical system. By converse transformation, we can obtain that:

\begin{Lemma}\label{L3.2}
Let $\varphi(t,\omega,z)$ be the random dynamical system generated by \eqref{RE-s}-\eqref{RE-f}. Then $T^{-1}(\theta_t\omega,\varphi(t,\omega,T(\omega,z))):=\tilde{\varphi}(t,\omega,z)$ is a random dynamical system. For any $z\in\mathbb{R}^{n+m}$,
the process $(t,\omega)\to\tilde{\varphi}(t,\omega,z)$ is a solution of \eqref{Equation-s}-\eqref{Equation-f}.
\end{Lemma}
\par
\noindent Hence, by a particular structure of transform $T$, if \eqref{RE-s}-\eqref{RE-f} has a stable or unstable foliation, so does \eqref{Equation-s}-\eqref{Equation-f}.
\par
 As we want to explore the relationship between the foliations and manifolds, we state the following results about the stable and unstable manifolds for \eqref{RE-s}-\eqref{RE-f}, similar to early works in \cite{DLS-2004,FLD}.

\begin{Lemma}(Random unstable manifold)\label{L3.3}
If the Lipschitz constant K, dichotomy parameters a, b satisfy the gap condition
$K(\frac{1}{\eta-b}+\frac{1}{a-\eta})<1$ with $b<\eta<a,$
then a Lipschitz  invariant random unstable manifold for the RDEs \eqref{RE-s}-\eqref{RE-f} exists, which is given by
\begin{equation}\label{manifold}
\mathcal{M}^u(\omega)=\{(\xi,
h^u(\xi,\omega)) |\; \xi \in \mathbb{R}^n \}
\end{equation}
where $h^u:\mathbb{R}^n\to\mathbb{R}^m$ is a Lipschitz continuous mapping that satisfies $h^u(0)=0$ and solves the following equation:
\begin{equation}\label{manifold-h}
h^u(\xi,\omega)=\int_{-\infty}^0e^{-Bs+\int_{s}^{0}z(\theta_\tau\omega)d\tau}G(\hat{x}(s,\omega; \xi),\hat{y}(s, \omega; \xi),\theta_s\omega)ds,\quad \hbox{for any}\quad \xi\in \mathbb{R}^n,
\end{equation}
\par

\noindent where $\hat{x}(t, \omega; \xi)$ and $\hat{y}(t,\omega; \xi)$ are the solutions of system
\eqref{RE-s}-\eqref{RE-f} with the following forms
$$
\left(
\begin{array}{l}
\hat{x}(t, \omega; \xi)\\
\hat{y}(t, \omega; \xi)
\end{array}
\right) =\left(
\begin{array}{l}
e^{At+\int_{0}^{t}z(\theta_\tau\omega)d\tau}\xi+\int_0^te^{A(t-s)+\int_{s}^{t}z(\theta_\tau\omega)d\tau}Fds\\
 \int_{-\infty}^t
e^{B(t-s)+\int_{s}^{t}z(\theta_\tau\omega)d\tau}Gds
\end{array}
\right)
$$
where $F=F(\hat{x}(s, \omega;\xi), \hat{y}(s, \omega; \xi), \theta_s\omega)$, $G=G(\hat{x}(s, \omega; \xi),\hat{y}(s, \omega; \xi),\theta_s\omega)$.  Furthermore,
$\tilde{\mathcal{M}}^u(\omega)=T^{-1}(\omega,\mathcal{M}^u(\omega))=\{(\xi,e^{z(\omega)}h^u(e^{-z(\omega)}\xi,\omega))|\xi\in\mathbb{R}^n\}$
is a Lipschitz unstable manifold of the stochastic differential system \eqref{Equation-s}-\eqref{Equation-f}.
\end{Lemma}
\par
Similar results on stable manifold can be obtained but we omit here.

\renewcommand{\theequation}{\thesection.\arabic{equation}}
\setcounter{equation}{0}

\section{Unstable foliation}

\quad\quad
To study system \eqref{RE-s}-\eqref{RE-f}, we define Banach spaces for a fixed $\eta$, $b<\eta<a$ as follows:
$$
C_{\eta}^{n,-}=\{\phi:(-\infty, 0]\to \mathbb{R}^n|\quad
\phi\;\hbox{is
 continous and}\; \sup\limits_{t\leq 0}e^{-\eta t-\int_{0}^{t}z(\theta_\tau\omega)d\tau}\left| \phi \right|<\infty \},
$$
$$
C_{\eta}^{n,+}=\{\phi:[0, +\infty)\to \mathbb{R}^n|\quad
\phi\;\hbox{is
 continous and}\; \sup\limits_{t\geq 0}e^{-\eta t-\int_{0}^{t}z(\theta_\tau\omega)d\tau}\left| \phi \right|<\infty \},
$$
with the norms
$$\|\phi\|_{C_{\eta}^{n,-}}=\sup\limits_{t\leq
0}e^{-\eta t-\int_{0}^{t}z(\theta_\tau\omega)d\tau} |\phi|,\quad
\hbox{and}\quad\|\phi\|_{C_{\eta}^{n,+}}=\sup\limits_{t\geq
0}e^{-\eta t-\int_{0}^{t}z(\theta_\tau\omega)d\tau} |\phi|,
$$
respectively. Analogously, we define Banach spaces $C_{\eta}^{m,-}$ and $C_{\eta}^{m,+}$ with the norms
$$\|\phi\|_{C_{\eta}^{m,-}}=\sup\limits_{t\leq
0}e^{-\eta t-\int_{0}^{t}z(\theta_\tau\omega)d\tau} |\phi|,\quad
\hbox{and}\quad\|\phi\|_{C_{\eta}^{m,+}}=\sup\limits_{t\geq
0}e^{-\eta t-\int_{0}^{t}z(\theta_\tau\omega)d\tau} |\phi|.\quad$$

 Let $C_\eta^\pm:=C_\eta^{n,\pm}\times C_\eta^{m,\pm}$, with norms $\|(x,
y)\|_{C_\eta^\pm}=\|x\|_{C_\eta^{n,\pm}}+\|y\|_{C_\eta^{m,\pm}}$, for
$(x, y)\in C_\eta^\pm$.
\par

Introduce
\begin{equation}\label{Foliation}
\mathcal{W}_\eta((\hat{x}_0, \hat{y}_0),\omega)=\{(\hat{x}_0^\ast,\hat{y}_0^\ast)\in \mathbb{R}^n\times \mathbb{R}^m| \;
\varphi(t,\omega, (\hat{x}_0,\hat{y}_0))-\varphi(t,\omega,(\hat{x}_0^\ast,\hat{y}_0^\ast))\in C_{\eta}^-\}.
\end{equation}
where $\varphi(t,\omega, (\hat{x}_0,\hat{y}_0))$ is the solution of the random system \eqref{RE-s}-\eqref{RE-f} as we denoted in Section 3. This is the set of all initial data through which the difference of two dynamical orbits are bounded by $e^{\eta t+\int_{0}^{t}z(\theta_\tau\omega)d\tau}.$
As we will prove later that $\mathcal{W}_\eta((\hat{x}_0, \hat{y}_0),\omega)$ is actually a fiber of the unstable foliation for the random
system \eqref{RE-s}-\eqref{RE-f}.
\par

Our main results about the existence of unstable foliation are present in the following.
\par
\begin{Theorem}(Unstable foliation)\label{Th4.1} Assume that the
Hypotheses A1-A2 hold. Take $\eta$ as the positive real number in the gap condition $K(\frac{1}{\eta-b}+\frac{1}{a-\eta})<1$. Then, the random dynamical system defined by \eqref{RE-s}-\eqref{RE-f} has a Lipschitz unstable foliation for which each unstable fiber can be represented as a graph
\begin{equation}\label{Foliation-mapping}
\mathcal{W}_\eta((\hat{x}_0, \hat{y}_0), \omega)=\{(\xi,
l(\xi, (\hat{x}_0, \hat{y}_0), \omega))| \; \xi\in\mathbb{R}^n\}.
\end{equation}
Here $(\hat{x}_0, \hat{y}_0)\in \mathbb{R}^n\times \mathbb{R}^m$, and the function
$l(\xi,(\hat{x}_0, \hat{y}_0), \omega)$ defined in \eqref{Foliation-l} is the graph mapping with Lipschitz constant satisfying
$$
Lip \; l \leq
\frac{K}{(\eta-b)\cdot(1-K(\frac{1}{\eta-b}+\frac{1}{a-\eta}))}.
$$
\end{Theorem}
\par
The proof of Theorem 4.1 based on the Lyapunov-Perron method will be presented after several useful Lemmas.
\par
Define the difference of two dynamical orbits of random system \eqref{RE-s}-\eqref{RE-f}
\begin{equation}\label{Equation-Psi}
\begin{array}{ll}
\phi(t)&=\varphi(t,\omega,(\hat{x}_0^\ast,\hat{y}_0^\ast))-\varphi(t,\omega,(\hat{x}_0,\hat{y}_0))\\
&=(\hat{x}(t,\omega,(\hat{x}_0^\ast,\hat{y}_0^\ast))-\hat{x}(t,\omega,(\hat{x}_0,\hat{y}_0)), \hat{y}(t,\omega,(\hat{x}_0^\ast,\hat{y}_0^\ast))-\hat{y}(t,\omega,
(\hat{x}_0,\hat{y}_0)))\\
&:=(u(t), v(t))
\end{array}
\end{equation}
with the initial condition
$$
\phi(0)=(u(0),v(0))=(\hat{x}_0^\ast-\hat{x}_0, \hat{y}_0^\ast-\hat{y}_0).
$$
Hence
$$
\begin{array}{l}
\hat{x}(t,\omega,(\hat{x}_0^\ast,\hat{y}_0^\ast))=u(t)+\hat{x}(t,\omega,(\hat{x}_0,\hat{y}_0)),\\
\hat{y}(t,\omega,(\hat{x}_0^\ast,\hat{y}_0^\ast))=v(t)+\hat{y}(t,\omega,(\hat{x}_0,\hat{y}_0)).
\end{array}
$$
Moreover, by using \eqref{RE-s}-\eqref{RE-f}, we find that $(u(t), v(t))$ satisfies the following equation:
\begin{eqnarray}
\frac{du}{dt}&=&Au+\Delta F(u,v, \theta_t\omega)+z(\theta_t\omega)u, \label{Equation-U}\\
\frac{dv}{dt}&=&Bv+\Delta G(u,v, \theta_t\omega)+z(\theta_t\omega)v, \label{Equation-V}
\end{eqnarray}
where
\begin{equation}\label{Equation-DF}
\begin{array}{lll}
\Delta F(u,v, \theta_t\omega)&=& F(u(t)+\hat{x}(t,\omega,(\hat{x}_0,\hat{y}_0)),v(t)+\hat{y}(t,\omega,(\hat{x}_0,\hat{y}_0)),\theta_t\omega)\\
&&-F(\hat{x}(t,\omega, (\hat{x}_0,\hat{y}_0)), \hat{y}(t,\omega,(\hat{x}_0,\hat{y}_0)), \theta_t\omega),
\end{array}
\end{equation}
\begin{equation}\label{Equation-DG}
\begin{array}{lll}
\Delta G(u,v, \theta_t\omega)&=& G(u(t)+\hat{x}(t,\omega,(\hat{x}_0,\hat{y}_0)),v(t)+\hat{y}(t,\omega,(\hat{x}_0,\hat{y}_0)),\theta_t\omega)\\
&&-G(\hat{x}(t,\omega, (\hat{x}_0,\hat{y}_0)), \hat{y}(t,\omega,(\hat{x}_0,\hat{y}_0)), \theta_t\omega),
\end{array}
\end{equation}
and   initial condition
$$
u(0)=u_0=\hat{x}_0^\ast-\hat{x}_0, ~
v(0)=v_0=\hat{y}_0^\ast-\hat{y}_0.
$$
Noted that the functions $\Delta F$ and $\Delta G$ also satisfy the Lipschitz condition with the same Lipschitz constant as $f$ or $g$.
\par
The following Lemma will offer the desired properties of the random function $\phi(t)=(u(t),v(t)).$
\par

\begin{Lemma}\label{L4.1} Suppose that $\phi(t)=(u(t), v(t))$ is in $C_{\eta}^-$. Then $\phi(t)$ is the solution of \eqref{Equation-U}-\eqref{Equation-V} with intial data $\phi(0)=(u_0,v_0)$  if and only if $\phi(t)$ satisfies :
\begin{equation}\label{Lemma 4.1-0}
\left(\begin{array}{l}u(t)\\v(t)\end{array}\right) = \left(
\begin{array}{c}
e^{At+\int_{0}^{t}z(\theta_\tau\omega)d\tau}u(0)+\int_0^te^{A(t-s)+\int_{s}^{t}z(\theta_\tau\omega)d\tau}\Delta F(u(s),v(s), \theta_s\omega)ds
\\
\int_{-\infty}^te^{B(t-s)+\int_{s}^{t}z(\theta_\tau\omega)d\tau}\Delta G(u(s), v(s), \theta_s\omega)ds
\end{array}
\right).
\end{equation}
\end{Lemma}
\par

\noindent {\bf Proof.} \emph{Necessity}. Suppose process $(u(t), v(t))$ solves system \eqref{Equation-U}-\eqref{Equation-V} with initial data $(u_0,v_0)$ and belong to Banach space $C_{\eta}^-$. Applying the variation of constants formula to system \eqref{Equation-U}-\eqref{Equation-V} for integral interval $r\leq t\leq0$,
$$u(t)=e^{A(t-r)+\int_{r}^{t}z(\theta_\tau\omega)d\tau}u(r)+\int_r^te^{A(t-s)+\int_{s}^{t}z(\theta_\tau\omega)d\tau}\Delta F(u(s),v(s), \theta_s\omega)ds,$$
$$v(t)=e^{B(t-r)+\int_{r}^{t}z(\theta_\tau\omega)d\tau}v(r)+\int_r^te^{B(t-s)+\int_{s}^{t}z(\theta_\tau\omega)d\tau}\Delta G(u(s),v(s), \theta_s\omega)ds.$$
We can check that the form of $u(t)$ is bounded under $\|\cdot\|_{C_{\eta}^{n,-}}$ by setting $r=0$.
\par
$$
\begin{array}{ll}
\|u(t)\|_{C_{\eta}^{n,-}}= \sup\limits_{t\leq0}e^{-\eta t-\int_{0}^{t}z(\theta_\tau\omega)d\tau}|u(t)|\\
~~~~\leq \sup\limits_{t\leq 0}\{e^{(a-\eta)t}|u(0)|+Ke^{-\eta t}\int_{t}^{0}e^{a(t-s)+\int_{s}^{0}z(\theta_\tau\omega)d\tau}(|u(s)|+|v(s)|)ds\}\\
~~~~\leq \sup\limits_{t\leq 0}\{e^{(a-\eta)t}|u(0)|+K\int_{t}^{0}e^{(a-\eta)(t-s)}(\|u(s)\|_{C_{\eta}^{n,-}}+\|v(s)\|_{C_{\eta}^{m,-}})ds\}\\
~~~~\leq |u(0)|+\frac{K}{a-\eta}(\|u(s)\|_{C_{\eta}^{n,-}}+\|v(s)\|_{C_{\eta}^{m,-}})<\infty
\end{array}
$$
\par
\noindent To make $(u(t),v(t))$ belong to the space $C_{\eta}^-$, apply the same ideas in the case of deterministic dynamical systems to find the appropriate form for $v(t)$, and notice that
\begin{eqnarray}
\|v(t)\|_{C_{\eta}^{m,-}}&=&\sup\limits_{t\leq0}e^{-\eta t}|e^{Bt}(e^{-Br+\int_{r}^{0}z(\theta_\tau\omega)d\tau}v(r)\notag\\
&&+\int_r^te^{-Bs+\int_{s}^{0}z(\theta_\tau\omega)d\tau}\Delta G(u(s),v(s), \theta_s\omega)ds)|\notag
\end{eqnarray}
Then for $t\leq0$ the following inequality holds and let $t\to-\infty$, we obtain
\begin{eqnarray}
|e^{-Br+\int_{r}^{0}z(\theta_\tau\omega)d\tau}v(r)+\int_r^te^{-Bs+\int_{s}^{0}z(\theta_\tau\omega)d\tau}\Delta G(u(s),v(s), \theta_s\omega)ds|\notag\\
\leq e^{(\eta-b)t}\|v(t)\|_{C_{\eta}^{m,-}}\to0\notag
\end{eqnarray}
which   implies for $t\leq0$,
$$v(r)=-e^{Br+\int_{0}^{r}z(\theta_\tau\omega)d\tau}\int_r^te^{-Bs+\int_{s}^{0}z(\theta_\tau\omega)d\tau}\Delta G(u(s),v(s), \theta_s\omega)ds$$
By taking limit for $t$, i.e., $t\to-\infty$ and replacing time variable r by t, we get
$$v(t)=\int_{-\infty}^te^{B(t-s)+\int_{s}^{t}z(\theta_\tau\omega)d\tau}\Delta G(u(s), v(s), \theta_s\omega)ds$$
Thus,   $(u(t), v(t))$ solving the system \eqref{Equation-U}-\eqref{Equation-V} in Banach space $C_{\eta}^-$ with initial data $(u_0,v_0)$ can be written as in \eqref{Lemma 4.1-0}.\par
\noindent \emph{Sufficiency}. By direct calculations, it is not hard to see that the process $\phi(t)=(u(t), v(t))$ is the solution of the system \eqref{Equation-U}-\eqref{Equation-V} if $\phi(t)$ can be written in form \eqref{Lemma 4.1-0} and is in $C_{\eta}^-$. This completes the proof of Lemma \ref{L4.1}.
\hfill{$\blacksquare$}
\par

From this Lemma,  we have the following Corollary.

\begin{Corollary}\label{C4.1} Assume that the Hypotheses A1-A2 hold.
Take $\eta$ as the positive real number. Then $(\hat{x}_0^\ast,\hat{y}_0^\ast)$ is in
$\mathcal{W}_\eta((\hat{x}_0, \hat{y}_0), \omega)$ if and only if
there exists a function $\phi(t)=(u(t),v(t))=(u(t, \omega, (\hat{x}_0, \hat{y}_0) ;u(0)), v(t, \omega, (\hat{x}_0, \hat{y}_0) ;u(0)))\in C_\eta^-$ satisfies
\eqref{Lemma 4.1-0}.
\end{Corollary}
\par

\begin{Lemma}\label{L4.2} Take $\eta>0$, $b<\eta<a$ so that they satisfy $K(\frac{1}{\eta-b}+\frac{1}{a-\eta})<1$. Given $u_0=\hat{x}_0^\ast-\hat{x}_0\in \mathbb{R}^n$, then the integral system \eqref{Lemma 4.1-0} has a
unique solution $\phi(\cdot)=\phi(\cdot,\omega, (\hat{x}_0, \hat{y}_0) ; u(0))$ in $C_\eta^-$ under the hypotheses A1-A2.
\end{Lemma}
\par
\noindent {\bf Proof.} To see this, for any $\phi(t)=(u(t), v(t))\in C_\eta^-$, introduce two operators
$\mathcal{J}_n: C_\eta^-\longrightarrow C_\eta^{n,-}$ and $\mathcal{J}_m: C_\eta^-\longrightarrow C_\eta^{m,-}$ by means of
$$
\mathcal{J}_n(\phi)[t]=e^{At+\int_{0}^{t}z(\theta_\tau\omega)d\tau}u(0)+\int_0^te^{A(t-s)+\int_{s}^{t}z(\theta_\tau\omega)d\tau}\Delta F(u(s),v(s), \theta_s\omega)ds,
$$
$$
\mathcal{J}_m(\phi)[t]
=\int_{-\infty}^te^{B(t-s)+\int_{s}^{t}z(\theta_\tau\omega)d\tau}\Delta G(u(s), v(s), \theta_s\omega)ds,~~~~~~~~~~
$$
for $t\leq0$ and define the mapping by\par
$$\mathcal{J}(\phi(\cdot)):=\binom{\mathcal{J}_n(\phi(\cdot))}{\mathcal{J}_m(\phi(\cdot))}.$$
It is easy to see that $\mathcal{J}$ is well-defined from $C_\eta^-$ into itself.
To this end, taking $\phi(t)=(u(t), v(t))\in C_\eta^-$, we have that
$$
\begin{array}{ll}
\|\mathcal{J}_n(\phi)[t]\|_{C_{\eta}^{n,-}}\\
~~\leq \sup\limits_{t\leq 0}\{e^{(a-\eta)t}|u(0)|+Ke^{-\eta t}\int_{t}^{0}e^{a(t-s)+\int_{s}^{0}z(\theta_\tau\omega)d\tau}(|u(s)|+|v(s)|)ds\}\\
~~\leq \sup\limits_{t\leq 0}\{e^{(a-\eta)t}|u(0)|+K\int_{t}^{0}e^{(a-\eta)(t-s)}(\|u(s)\|_{C_{\eta}^{n,-}}+\|v(s)\|_{C_{\eta}^{m,-}})ds\}\\
~~\leq |u(0)|+\frac{K}{a-\eta}(\|u(s)\|_{C_{\eta}^{n,-}}+\|v(s)\|_{C_{\eta}^{m,-}})\\
~~=|u(0)|+\frac{K}{a-\eta}\|\phi(t)\|_{C_\eta^-}
\end{array}
$$
\noindent and
$$
\begin{array}{ll}
\|\mathcal{J}_m(\phi)[t]\|_{C_{\eta}^{m,-}}
&\leq \sup\limits_{t\leq 0}\{Ke^{-\eta t}\int_{-\infty}^{t}e^{b(t-s)+\int_{s}^{0}z(\theta_\tau\omega)d\tau}(|u(s)|+|v(s)|)ds\}\\
&\leq \sup\limits_{t\leq 0}\{K\int_{-\infty}^{t}e^{(b-\eta)(t-s)}\|\phi(s)\|_{C_\eta^-}\}\\
&=\frac{K}{\eta-b}\|\phi(s)\|_{C_\eta^-}.\\
\end{array}
$$
\noindent Hence, by the definition of $\mathcal{J}$, we obtain
$$\mathcal{J}(\phi(t))\leq|u(0)|+(\frac{K}{a-\eta}+\frac{K}{\eta-b})\|\phi(t)\|_{C_\eta^-}.$$
\par
\noindent Thus, we conclude that $\mathcal{J}$ maps $C_\eta^-$ into itself.\par
\noindent Further, we will show that the mapping $\mathcal{J}$ is contractive. To see this, taking any $\phi=(u, v)\in C_\eta^-$ and
$\hat{\phi}=(\hat{u},\hat{v})\in C_\eta^-$, then
\begin{equation}\label{Lemma 4.2-1}
\begin{array}{ll}
&\|\mathcal{J}_n(\phi)- \mathcal{J}_n(\hat{\phi})\|_{C_\eta^{n,-}}\\
=&\|\int_{0}^t e^{A(t-s)+\int_{s}^{t}z(\theta_\tau\omega)d\tau}[\Delta F(u(s),
v(s),\theta_s\omega)- \Delta
F(\hat{u}(s),\hat{v}(s),\theta_s\omega)]ds\|_{C_\eta^{n,-}}\\
=&\|\int_{0}^te^{A(t-s)+\int_{s}^{t}z(\theta_\tau\omega)d\tau}[F(u(s)+\hat{x}(s,\omega, (\hat{x}_0,\hat{y}_0)), v(s)+\hat{y}(s,\omega, (\hat{x}_0,\hat{y}_0)),\theta_s\omega)\\
&~~~~~- F(\hat{u}(s)+\hat{x}(s,\omega, (\hat{x}_0,\hat{y}_0)), \hat{v}(s)+\hat{y}(s,\omega, (\hat{x}_0,\hat{y}_0)),\theta_s\omega)]ds
\|_{C_\eta^{n,-}}\\
\leq &\sup\limits_{t\leq 0}\{Ke^{-\eta t}\int_{t}^{0}e^{a(t-s)+\int_{s}^{0}z(\theta_\tau\omega)d\tau}
(|u(s)-\hat{u}(s)|+|v(s)-\hat{v}(s)|)ds\}\\
\leq&\sup\limits_{t\leq 0}\{K\int_{t}^{0}e^{(a-\eta)(t-s)}\|\phi-\hat{\phi}\|_{C_\eta^-}ds\}\\
\leq&\frac{K}{a-\eta}\|\phi-\hat{\phi}\|_{C_\eta^-},
\end{array}
\end{equation}
and
\begin{equation}\label{Lemma 4.2-2}
\begin{array}{ll}
&\|\mathcal{J}_m(\phi)- \mathcal{J}_m(\hat{\phi})\|_{C_\eta^{m,-}}\\
=&\|\int_{-\infty}^te^{B(t-s)+\int_{s}^{t}z(\theta_\tau\omega)d\tau}[\Delta
G(u(s), v(s),\theta_s\omega)-\Delta G(\hat{u}(s),\hat{v}(s),\theta_s\omega)]ds\|_{C_\eta^{m,-}}\\
=&\|\int_{-\infty}^te^{B(t-s)+\int_{s}^{t}z(\theta_\tau\omega)d\tau}[G(u(s)+\hat{x}(s,\omega, (\hat{x}_0, \hat{y}_0), v(s)+\hat{y}(s,\omega, (\hat{x}_0, \hat{y}_0)),\theta_s\omega)\\
&~~~~~- G(\hat{u}(s)+\hat{x}(s,\omega, (\hat{x}_0, \hat{y}_0)), \hat{v}(s)+\hat{y}(s,\omega, (\hat{x}_0, \hat{y}_0)),\theta_s\omega)]ds\|_{C_\eta^{m,-}}\\
\leq &\sup\limits_{t\leq 0}Ke^{-\eta t}\int_{-\infty}^{t}e^{b(t-s)+\int_{s}^{0}z(\theta_\tau\omega)d\tau}
(|u(s)-\hat{u}(s)|+|v(s)-\hat{v}(s)|)ds\\
\leq &\frac{K}{\eta-b}\|\phi-\hat{\phi}\|_{C_\eta^-}.
\end{array}
\end{equation}
\noindent Hence, by \eqref{Lemma 4.2-1} and \eqref{Lemma 4.2-2}
\begin{equation}\label{Lemma 4.2-3}
\begin{array}{ll}
\|\mathcal{J}(\phi)-\mathcal{J}(\hat{\phi})\|_{C_\eta^-}
&=\|\mathcal{J}_n(\phi)-\mathcal{J}_n(\hat{\phi})\|_{C_\eta^{n,-}}+\|\mathcal{J}_m(\phi)-\mathcal{J}_m(\hat{\phi})\|_{C_\eta^{m,-}}\\
&\leq (\frac{K}{a-\eta}+\frac{K}{\eta-b})\|\phi-\hat{\phi}\|_{C_\eta^-}.
\end{array}
\end{equation}
Put the constant
\begin{equation}\label{Lemma 4.2-4}
\rho(a,b,K)=\frac{K}{a-\eta}+\frac{K}{\eta-b}.
\end{equation}
Then
\begin{equation}\label{Lemma 3.2-5}
\|\mathcal{J}(\phi)-\mathcal{J}(\hat{\phi})\|_{C_\eta^-}\leq
\rho(a,b,K)\|\phi-\hat{\phi}\|_{C_\eta^-}.
\end{equation}
\par
\noindent By the assumption,
$0<\rho(a,b,K)<1.$
Hence the map $\mathcal{J}(\phi)$ is contractive in $C_\eta^-$ uniformly with respect to $(\omega, (\hat{x}_0, \hat{y}_0); u(0))$. By the uniform contraction mapping principle,
we have that the mapping
$\mathcal{J}(\phi)=\mathcal{J}(\phi,\omega, (\hat{x}_0, \hat{y}_0); u(0))$ has a unique fixed point for each $u(0)\in \mathbb{R}^n$,
which still denoted by $$\phi(\cdot)
=\phi(\cdot, \omega, (\hat{x}_0, \hat{y}_0); u(0)) \in C_\eta^-.$$ Namely, $\phi(\cdot, \omega, (\hat{x}_0, \hat{y}_0); u(0))\in C_\eta^-$ is a unique solution of the system \eqref{Lemma 4.1-0} with the initial data $u(0)$.\hfill{$\blacksquare$}
\par


Lemma \ref{L4.2} ensures the existence and uniqueness of solution of the system \eqref{Lemma 4.1-0} for each given initial value. In fact, following Lemma indicates that the solution of the system \eqref{Lemma
4.1-0}, i.e., $\phi(t)=\phi(t, \omega, (\hat{x}_0, \hat{y}_0); u(0))$ has continuous dependence on initial conditions.
\par
\begin{Lemma}\label{L4.3}
Assume the same conditions as stated in Lemma \ref{L4.2}. Let $\phi(t)=\phi(t, \omega, (\hat{x}_0, \hat{y}_0); u(0))$ be the unique solution of the system \eqref{Lemma 4.1-0} in $C_{\eta}^-$. Then for every $u(0)$ and $\tilde{u}(0)$ in $\mathbb{R}^n$, we have
\begin{equation}\label{Lemma 4.3.1}
\|\phi(t, \omega, (\hat{x}_0, \hat{y}_0); u(0))-\phi(t, \omega, (\hat{x}_0, \hat{y}_0); \tilde{u}(0))\|_{C_{\eta}^-}
\leq \frac{1}{1-\rho(a,b,K)}|u(0)-\tilde{u}(0)|,
\end{equation}
where $\rho(a,b,K)$ is defined as
\eqref{Lemma 4.2-4}.
\end{Lemma}
\par
\noindent{\bf Proof.} Taking any $u(0)$ and $\tilde{u}(0)$ in $\mathbb{R}^n$, we write $u(t,\omega;u(0))$ instead of $u(t,\omega,(\hat{x}_0, \hat{y}_0);u(0))$ and $v(t,\omega;u(0))$ instead of $v(t,\omega,(\hat{x}_0, \hat{y}_0);u(0))$ for simplicity. We have for the fixed point $\phi$ the estimate:
$$
\begin{array}{ll}
&\|\phi(t, \omega, (\hat{x}_0, \hat{y}_0);u(0))-\phi(t, \omega, (\hat{x}_0, \hat{y}_0); \tilde{u}(0))\|_{C_{\eta}^-}\\
=&|u(t,\omega;u(0))-u(t,\omega;\tilde u(0))|+|v(t,\omega;u(0))-v(t,\omega;\tilde u(0))|\\
\leq &
|u(0)-\tilde{u}(0)|+\rho(a,b,K)\|\phi(t, \omega, (\hat{x}_0, \hat{y}_0);u(0))-\phi(t, \omega, (\hat{x}_0, \hat{y}_0); \tilde{u}(0))\|_{C_\eta^-}.
\end{array}
$$
\noindent Then we obtain \eqref{Lemma 4.3.1} by transposition.\hfill{$\blacksquare$}

\par
In the following, for every $\xi\in \mathbb{R}^n$, we define a function:
\begin{equation}\label{Foliation-l}
 l(\xi,(\hat{x}_0, \hat{y}_0),\omega):=\hat{y}_0+\int_{-\infty}^0e^{-Bs+\int_{s}^{0}z(\theta_\tau\omega)d\tau}\Delta G(u,v, \theta_s\omega)ds.
\end{equation}
with $u=u(s, \omega, (\hat{x}_0, \hat{y}_0) ; (\xi-\hat{x}_0))$, $v=v(s, \omega, (\hat{x}_0, \hat{y}_0);(\xi-\hat{x}_0))$.\par
\par
\noindent{\bf Proof of Theorem \ref{Th4.1}.}\par
\noindent From \eqref{Lemma 4.1-0}, we deduce that
$$
\left(
\begin{array}{l}
\hat{x}_0^\ast-\hat{x}_0\\
\hat{y}_0^\ast-\hat{y}_0
\end{array}
\right ) = \left(
\begin{array}{c}
\hat{x}_0^\ast-\hat{x}_0\\
\int_{-\infty}^0e^{-Bs+\int_{s}^{0}z(\theta_\tau\omega)d\tau}
\Delta G(u(s),v(s),\theta_s\omega)ds
\end{array}
\right ).
$$
As a sequence,
\begin{eqnarray}
\hat{y}_0^\ast&=&\hat{y}_0+\int_{-\infty}^0\exp(-Bs+\int_{s}^{0}z(\theta_\tau\omega)d\tau)\cdot\notag\\
&&~~~~~~~~~~~~~\Delta G(u(s, \omega, (\hat{x}_0, \hat{y}_0); u(0)),v(s,\omega, (\hat{x}_0, \hat{y}_0); u(0)), \theta_s\omega)ds\notag\\
&=&\hat{y}_0+\int_{-\infty}^0\exp(-Bs+\int_{s}^{0}z(\theta_\tau\omega)d\tau)\cdot\notag\\
&&~~~~~~~~~~~~\Delta G(u(s, \omega, (\hat{x}_0, \hat{y}_0);\hat{x}_0^\ast-\hat{x}_0 ),v(s, \omega, (\hat{x}_0, \hat{y}_0) ; (\hat{x}_0^\ast-\hat{x}_0)), \theta_s\omega)ds,\notag
\end{eqnarray}
We find that above function just is $l(\xi,(\hat{x}_0, \hat{y}_0),\omega)$ if we take
$\hat{x}_0^\ast$ as $\xi$ in $\mathbb{R}^n$. Then according to Corollary \ref{C4.1}, Lemma
\ref{L4.2}, \eqref{Foliation} and \eqref{Foliation-l}, we see that
$$
\mathcal{W}_\eta((\hat{x}_0, \hat{y}_0), \omega)=\{(\xi,
l(\xi,(\hat{x}_0, \hat{y}_0),\omega))| \; \xi\in \mathbb{R}^n\},
$$
which immediately shows a fiber of the unstable foliation can be represented as graph of a function.\par\noindent In addition, for any $\xi$ and $\tilde{\xi}$ in $\mathbb{R}^n$, via
\eqref{Foliation-l} and Lemma \ref{L4.3},
$$
\begin{array}{ll}
&|l(\xi,(\hat{x}_0, \hat{y}_0),\omega)-l(\tilde{\xi},(\hat{x}_0, \hat{y}_0),\omega)|\\
\leq&\frac{K}{\eta-b}\|\phi(\cdot,\omega,(\hat{x}_0, \hat{y}_0); \xi-\hat{x}_0)-\phi(\cdot,\omega, (\hat{x}_0, \hat{y}_0);\tilde{\xi}-\hat{x}_0)\|_{C_\eta^{-}}\\
\leq&
\frac{K}{\eta-b}\cdot\frac{1}{1-\rho(a,b,K)}|\xi-\tilde{\xi}|.
\end{array}
$$
This shows that $l(\xi,(\hat{x}_0, \hat{y}_0),\omega)$ is Lipschitz continuous with respect to variable $\xi$.
\par \noindent Thus  the proof is complete.\hfill{$\blacksquare$}
\par


\begin{Remark}\label{R4.1}
Notice that the relationship between the solutions of system \eqref{Equation-s}-\eqref{Equation-f} and \eqref{RE-s}-\eqref{RE-f}, the original stochastic system also has an unstable foliation under the conditions of Theorem \ref{Th4.1}, and every unstable fiber is represented as
\begin{eqnarray}
\tilde{\mathcal{W}}_\eta((x_0, y_0), \omega)&=&T^{-1}(\omega,\mathcal{W}_\eta((\hat{x}_0, \hat{y}_0), \omega))\notag\\
&=&(\{(\xi,e^{z(\omega)}l(e^{-z(\omega)}\xi,(x_{0}e^{-z(\omega)},y_0e^{-z(\omega)}),\omega))|\xi\in\mathbb{R}^n\}.\notag
\end{eqnarray}
Different from the case of Brownian noise, the dynamical orbits in $\tilde{\mathcal{W}}_\eta((x_0, y_0), \omega)$ are c\'adl\'ag and adapted.
\end{Remark}\par
In what follows we are going to prove that if dynamical orbits of \eqref{RE-s}-\eqref{RE-f} start from the
same unstable fiber, then they will approach each other exponentially in backward time.  \par

\begin{Theorem}(Properties of unstable foliation)\label{Th4.2} Assume that the
Hypotheses A1-A2 hold. Take $\eta>0$, $b<\eta<a$  so that they satisfy the gap condition  $K(\frac{1}{\eta-b}+\frac{1}{a-\eta})<1$. Then, the Lipschitz unstable foliation for \eqref{RE-s}-\eqref{RE-f} obtained in Theorem \ref{Th4.1} has the following properties:
\par
\noindent(i)\quad 
the dynamical orbits which start from the
same fiber are exponentially
approaching each other in backward time. In other words, for every two points $(\hat{x}_0^1,
\hat{y}_0^1)$ and $(\hat{x}_0^2, \hat{y}_0^2)$
in a same fiber $\mathcal{W}_\eta((\hat{x}_0, \hat{y}_0),\omega)$,
\begin{equation}\label{Foliation-attracting}
\begin{array}{ll}
|\varphi(t,\omega, (\hat{x}_0^1,
\hat{y}_0^1))-\varphi(t,\omega, (\hat{x}_0^2,
\hat{y}_0^2))|&\leq \frac{e^{\eta
t+\int_{0}^{t}z(\theta_\tau\omega)d\tau}}{1-\rho(a,b,K)}\cdot
|\hat{x}_0^1-\hat{x}_0^2|\\
& =O(e^{\eta t}),~~~~~~~~~~ \forall~~ \omega, ~~ as \quad t\to-\infty.
\end{array}
\end{equation}
\noindent(ii)\quad its unstable fiber is invariant in the sense of cocycle, i.e.,
$$
\varphi(t, \omega, \mathcal{W}_\eta((\hat{x}_0, \hat{y}_0), \omega))\subset
\mathcal{W}_\eta(\varphi(t,\omega, (\hat{x}_0, \hat{y}_0)), \theta_t\omega).
$$
\end{Theorem}

\noindent{\bf Proof.} (i) In view of Corollary \ref{C4.1} as well as the same argument in the proof of Lemma \ref{L4.2},
we find that
\begin{equation}\label{Theorem 4.1-1}
\begin{array}{lll}
\|\phi(\cdot)\|_{C_\eta^-}&=&\|u(\cdot)\|_{C_\eta^{n,-}}+\|v(\cdot)\|_{C_\eta^{m,-}}\\
&\leq&\|e^{At+\int_{0}^{t}z(\theta_\tau\omega)d\tau}u(0)\|_{C_\eta^{n,-}}\\
&&+\|\int_0^te^{A(t-s)+\int_{s}^{t}z(\theta_\tau\omega)d\tau}\Delta F(u(s), v(s),
\theta_s\omega)ds\|_{C_\eta^{n,-}}\\&&+\|\int_{-\infty}^t
e^{B(t-s)+\int_{s}^{t}z(\theta_\tau\omega)d\tau}
\Delta
G(u(s), v(s), \theta_s\omega)ds\|_{C_\eta^{m,-}}\\&\leq & |u(0)|
+\frac{K}{a-\eta}\|\phi(\cdot)\|_{C_\eta^-}
+\frac{K}{\eta-b}\|\phi(\cdot)\|_{C_\eta^-}\\
&\leq & |u(0)|+\rho(a,b,K)\|\phi(\cdot)\|_{C_\eta^-},
\end{array}
\end{equation}
where $\phi$ is defined as \eqref{Equation-Psi}. Then it follows from
\eqref{Theorem 4.1-1} that
$$
\|\phi(\cdot)\|_{C_\eta^-}\leq
\frac{1}{1-\rho(a,b,K)}|u(0)|,
$$
which implies immediately that
\begin{equation}\label{Theorem 4.2-2}
|\varphi(t,\omega,
(\hat{x}_0^\ast,\hat{y}_0^\ast))-\varphi(t,\omega,(\hat{x}_0,\hat{y}_0))|\leq \frac{e^{\eta
t+\int_{0}^{t}z(\theta_\tau\omega)d\tau}}{1-\rho(a,b,K)}\cdot
|u(0)|, \quad \forall \quad t\leq 0.
\end{equation}
\par
\noindent Hence, for every set of two points $(\hat{x}_0^1, \hat{y}_0^1)$ and
$(\hat{x}_0^2, \hat{y}_0^2)$ from the same fiber
$\mathcal{W}_\eta((\hat{x}_0,\hat{y}_0), \omega)$,  as both of them satisfy \eqref{Theorem 4.2-2}, we have
$$
|\varphi(t,\omega, (\hat{x}_0^1,
\hat{y}_0^1))-\varphi(t,\omega,
(\hat{x}_0,\hat{y}_0))|\leq \frac{e^{\eta
t+\int_{0}^{t}z(\theta_\tau\omega)d\tau}}{1-\rho(a,b,K)}\cdot
|u(0)|, \quad \forall \quad t\leq 0,
$$
and
$$
|\varphi(t,\omega, (\hat{x}_0^2,
\hat{y}_0^2))-\varphi(t,\omega,
(\hat{x}_0,\hat{y}_0))|\leq \frac{e^{\eta
t+\int_{0}^{t}z(\theta_\tau\omega)d\tau}}{1-\rho(a,b,K)}\cdot
|u(0)|, \quad \forall \quad t\leq 0.
$$
These imply that \eqref{Foliation-attracting} holds apparently.\par Notice that
$$\lim_{t \to \pm\infty}\frac{1}{t}\int_{0}^tz(\theta_\tau\omega)d\tau=0, \quad\omega \in \Omega$$}In other words, $\int_{0}^tz(\theta_\tau\omega)d\tau$ has a sublinear growth rate which is increasing slowly than linear increasing, thus, $e^{\int_{0}^{t}z(\theta_\tau\omega)d\tau}$ does not change the exponential convergence of solutions starting at the same fiber, the proof of (i) is complete.\par
\par
\noindent(ii) To prove the fiber   invariance, taking a fiber
$\mathcal{W}_\eta((\hat{x}_0,\hat{y}_0), \omega)$ arbitrarily, we need to show
that $$\varphi(\tau, \omega, \mathcal{W}_\eta((\hat{x}_0,\hat{y}_0), \omega))\subset
\mathcal{W}_\eta(\varphi(\tau, \omega, ( \hat{x}_0,\hat{y}_0)), \theta_\tau\omega).$$ Let $(\hat{x}_0^\ast,
\hat{y}_0^\ast)\in \mathcal{W}_\eta((\hat{x}_0,\hat{y}_0),
\omega)$. We have $\varphi(\cdot, \omega, (\hat{x}_0^\ast,
\hat{y}_0^\ast))-\varphi(\cdot, \omega, (\hat{x}_0,\hat{y}_0))\in
C_{\eta}^-$ from \eqref{Foliation}, which implies that
$$\varphi(\cdot+\tau,
\omega, (\hat{x}_0^\ast,
\hat{y}_0^\ast))-\varphi(\cdot+\tau, \omega, (\hat{x}_0,\hat{y}_0))\in C_{\eta}^-.
$$
Thus according to the cocycle property
$$
\begin{array}{l}
\varphi(\cdot +\tau, \omega, (\hat{x}_0^\ast,
\hat{y}_0^\ast))=\varphi(\cdot, \theta_\tau\omega,
\varphi(\tau, \omega, (\hat{x}_0^\ast,
\hat{y}_0^\ast))),\\
\varphi(\cdot +\tau, \omega, (\hat{x}_0,\hat{y}_0))=\varphi(\cdot, \theta_\tau\omega,
\varphi(\tau, \omega, (\hat{x}_0,\hat{y}_0))),
\end{array}
$$
hence
$\varphi(\cdot, \theta_\tau\omega,
\varphi(\tau, \omega, (\hat{x}_0^\ast,
\hat{y}_0^\ast)))-\varphi(\cdot, \theta_\tau\omega,
\varphi(\tau, \omega, (\hat{x}_0,\hat{y}_0)))\in C_{\eta}^-$.\par
\noindent Then we have $\varphi(\tau, \omega, (\hat{x}_0^\ast,
\hat{y}_0^\ast))\in
\mathcal{W}_\eta(\varphi(\tau, \omega, (\hat{x}_0,\hat{y}_0)), \theta_\tau\omega)$. The proof is complete.
\hfill{$\blacksquare$}
\par

\begin{Remark}\label{R4.2} Under the same conditions presented in Theorem \ref{Th4.2}, the unstable foliation of original stochastic system \eqref{Equation-s}-\eqref{Equation-f} is also invariant because of the nature of the random transformation $T$. Furthermore, note that $t\to z(\theta_t\omega)$ has a sublinear growth rate guaranteed by Lemma \ref{L3.1}. Thus, the transform $T^{-1}(\theta_t\omega,\cdot)$ does not change the exponential convergence of dynamical orbits of the system \eqref{Equation-s}-\eqref{Equation-f} in backward time starting from the same fiber.
\end{Remark}


\begin{Remark}\label{R4.3} In addition, by early works as well as the results of this paper, we   see that the unstable foliation and unstable manifold are the useful tools describing different aspects of the dynamics for stochastic systems with multiplicative non-Gaussian noise.
\end{Remark}
\par

\begin{Remark}\label{R4.4} Usual gap condition $\frac{K}{a-\eta}+\frac{K}{\eta-b}<1$ given in \cite{CLL} and \cite{DLS-2004} only indicates the existence of the mapping $l$ of the
unstable foliation for the random system
\eqref{RE-s}-\eqref{RE-f}. To ensure dynamical orbits starting from the same
fiber exponentially approaching each other in backward
time, we require $\eta>0$ additionally.
\par
More precisely, for the existence of unstable foliation, we need: (i) $a-\eta>0$ in \eqref{Lemma 4.2-1} and
$\frac{K}{a-\eta}+\frac{K}{\eta-b}<1$ in \eqref{Lemma
4.2-4}; for the exponentially approaching of dynamical orbits, need: (ii)
$\eta>0$. Therefore, a simple choice is that $\eta=\frac{1}{p}a$
with $p>1$. By directly calculation, we find that the corresponding exponentially approaching
rate is less if we choose the gap
condition is larger. So it is unfortunate that we can not obtain
the optimal gap condition with the optimal rate $\eta$. For simplicity, we only require $\eta>0$ and dose not specify the exact value that $\eta$ takes.
\end{Remark}
\par
The link between the unstable foliation and unstable manifold is presented in the following theorem (refer to \cite{Cdz} for the case of additive Brownian noise.)\par

\begin{Theorem}(Geometric structures of the unstable foliation)\label{Th4.3} Assume that the Hypotheses A1-A2 hold. Take $\eta>0$, $b<\eta<a$. Let $\mathcal{M}^u(\omega)$
 and $\mathcal{W}_\eta((\hat{x}_0,
\hat{y}_0), \omega)$ be the unstable manifold and a fiber
of the unstable foliation for the random system
\eqref{RE-s}-\eqref{RE-f}, which are well defined in
\eqref{manifold} and \eqref{Foliation-mapping}, respectively.
Put
$$
\overset{\text{p}}{\mathcal{W}_\eta}(\omega):=\{\mathcal{W}_\eta((\hat{x}_0,
\hat{y}_0), \omega)\;|\; \hat{y}_0-h^u(\hat{x}_0, \omega):=p \in \mathbb{R}^m,\;
(\hat{x}_0,
\hat{y}_0)\in \mathbb{R}^n\times \mathbb{R}^m\},
$$
where $h^u(\hat{x}_0, \omega)$ is defined in \eqref{manifold-h}.
%
%
Then
\par
\noindent(i)\quad if $p=0$,
$\overset{\text{p}}{\mathcal{W}_\eta}(\omega)$ is just
the unstable manifold.
\par
\noindent(ii)\quad for any $p, q\in \mathbb{R}^m$ and $p\neq q$, the unstable fiber
$\overset{\text{p}}{\mathcal{W}_\eta}(\omega)$
parallels to the unstable fiber
$\overset{\text{q}}{\mathcal{W}_\eta}(\omega)$.
\par
\noindent(iii)\quad the unstable fiber $\mathcal{W}_\eta((\hat{x}_0, \hat{y}_0), \omega)$ is
just the unstable manifold $\mathcal{M}^u(\omega)$ if we choose an arbitrarily point $(\hat{x}_0, \hat{y}_0)$ from the unstable
foliation and this chosen point also belongs to unstable manifold $\mathcal{M}^u(\omega)$.
\par
\noindent(iv)\quad the unstable fiber $\mathcal{W}_\eta((\hat{x}_0, \hat{y}_0), \omega)$ and unstable manifold $\mathcal{M}^u(\omega)$ are parallel if arbitrarily taken point $(X_0, Y_0)$ of the unstable
foliation is not in the unstable manifold $\mathcal{M}^u(\omega).$
\end{Theorem}
\par

\noindent{\bf Proof .}\quad
It follows from \eqref{Lemma 4.1-0} that, for any
$(\hat{x}_0^\ast,\hat{y}_0^\ast)\in \mathbb{R}^n\times\mathbb{R}^m$, we have
$$
\hat{y}_0^\ast-\hat{y}_0 =
\int_{-\infty}^0e^{-Bs+\int_{s}^{0}z(\theta_\tau\omega)d\tau}\Delta
G(u(s), v(s), \theta_s\omega)ds,
$$
which suggests that
\begin{equation}\label{Theorem 4.3-1-1}
\begin{array}{ll}
&\hat{y}_0^\ast-
\int_{-\infty}^0e^{-Bs+\int_{s}^{0}z(\theta_\tau\omega)d\tau}G(\hat{x}(s,
\omega; \hat{x}_0^\ast),  \hat{y}(s, \omega;
\hat{x}_0^\ast), \theta_s\omega)ds\\
=&\hat{y}_0 -
\int_{-\infty}^0e^{-Bs+\int_{s}^{0}z(\theta_\tau\omega)d\tau}G(\hat{x}(s,\omega; \hat{x}_0),  \hat{y}(s, \omega; \hat{x}_0),\theta_s\omega)ds.
\end{array}
\end{equation}
Namely,
\begin{equation}\label{Theorem 4.3-1-2}
\hat{y}_0^\ast-h^u(\hat{x}_0^\ast,\omega)=\hat{y}_0-h^u(\hat{x}_0,\omega),
\end{equation}
where $h^u(\cdot, \omega)$ is defined as
\eqref{manifold-h}.
\par
\noindent If we take an arbitrary point $(\hat{x}_0, \hat{y}_0)$ from the unstable foliation,
then there exists $p\in \mathbb{R}^m$ such that
$$
\hat{y}_0-h^u(\hat{x}_0,\omega)=p.
$$
\par
\noindent When $p=0$, then $(\hat{x}_0, \hat{y}_0)$ belong to the unstable manifold
$\mathcal{M}^u(\omega)$, we obtain from \eqref{Theorem
4.3-1-2}
$$
\hat{y}_0^\ast-h^u(\hat{x}_0^\ast,\omega)=0, \quad
\hbox{for any}\quad  \hat{x}_0^\ast\in \mathbb{R}^n.
$$
Thus,
$\overset{\text{0}}{\mathcal{W}_\eta}(\omega)=\mathcal{M}^u(\omega)$.
\par
\noindent When $p\neq 0$, then $(\hat{x}_0, \hat{y}_0)$ is not in the unstable manifold
$\mathcal{M}^u(\omega)$. Then it immediately yields from
\eqref{Theorem 4.3-1-2} that
$$
\hat{y}_0^\ast-h^u(\hat{x}_0^\ast,\omega)=p\neq0,\quad
\hbox{for any}\quad  \hat{x}_0^\ast\in \mathbb{R}^n.
$$
Thus $(\hat{x}_0^\ast, \hat{y}_0^\ast)$ falls into
$\overset{\text{p}}{\mathcal{W}_\eta}(\omega)$ that
parallels  to the unstable manifold
$\mathcal{M}(\omega)=\overset{\text{0}}{\mathcal{W}_\eta}(\omega)$.
And apparently, for $p, q\in \mathbb{R}^m$ and $p\neq q$, the
$\overset{\text{p}}{\mathcal{W}_\eta}(\omega)$
parallels to
$\overset{\text{q}}{\mathcal{W}_\eta}(\omega)$. The
proof is completed. $\blacksquare$\par

\begin{Remark}\label{R4.5} From Theorem \ref{Th4.3}, we have a clear idea of geometric structure of the unstable foliation: (i) fibers of the unstable foliation parallel to each other; (ii) the unstable manifold is a special fiber. Namely, if we take an arbitrary point from the a fiber and this chosen point just falls in the unstable manifold, then the fiber just be the unstable manifold itself. Finally, what needs to explain is that fiber paralleling with each other here means
that the two fibers have parallel tangent lines at each corresponding horizontal point.
\end{Remark}
\par

Analogously, we  also obtain the corresponding results on the stable foliation stated in the following theorem without proof.
\par

\begin{Theorem}(stable foliation)\label{Th4.4} Assume that the
Hypotheses A1-A2 hold. Take $\gamma<0$, $b<\gamma<a$ so that they satisfy the gap condition $K(\frac{1}{\gamma-b}+\frac{1}{a-\gamma})<1$. Then, \par
\noindent(i) the random dynamical system defined by \eqref{RE-s}-\eqref{RE-f} have an invariant Lipschitz stable foliation for which every fiber is represented as a graph
\begin{equation}\label{Foliation-mapping-stable}
\mathcal{W}_\gamma((\hat{x}_0, \hat{y}_0), \omega)=\{(l(\zeta,(\hat{x}_0, \hat{y}_0), \omega),\zeta)| \; \zeta\in\mathbb{R}^m\},
\end{equation}
where $(\hat{x}_0, \hat{y}_0)\in \mathbb{R}^n\times \mathbb{R}^m$.  The function
$l(\zeta, (\hat{x}_0, \hat{y}_0), \omega)$ is the graph mapping with Lipschitz constant satisfying
$$
Lip \; l \leq
\frac{K}{(a-\gamma)\cdot(1-K(\frac{1}{\gamma-b}+\frac{1}{a-\gamma}))},
$$
where $l(\zeta, (\hat{x}_0, \hat{y}_0), \omega)$ is defined as\par
$\begin{array}{ll}
l(\zeta,(\hat{x}_0, \hat{y}_0),\omega):=\hat{x}_0
+\int_\infty^0 \exp(-As+\int_{s}^{0}z(\theta_\tau\omega)d\tau)\cdot\\
~~~~~~~~~~~~~\Delta F(u(s, \omega, (\hat{x}_0, \hat{y}_0) ; (\zeta-\hat{y}_0)),v(s, \omega, (\hat{x}_0, \hat{y}_0);
(\zeta-\hat{y}_0)), \theta_s\omega)ds.
\end{array}$\par
\noindent Furthermore, by an inverse transformation,
\begin{eqnarray}
\tilde{\mathcal{W}}_\gamma((x_0, y_0), \omega)&=&T^{-1}(\omega,\mathcal{W}_\gamma((\hat{x}_0, \hat{y}_0), \omega))\notag\\
&=&\{(e^{z(\omega)}l(e^{-z(\omega)}\zeta,(x_{0}e^{-z(\omega)},y_0e^{-z(\omega)}),\omega),\zeta)|\zeta\in\mathbb{R}^m\}\notag
\end{eqnarray}
is a Lipschitz stable foliation of the original  stochastic system \eqref{Equation-s}-\eqref{Equation-f}.\par
\noindent(ii) the Lipschitz stable foliation for \eqref{RE-s}-\eqref{RE-f} obtained above have the following property:
the dynamical orbits which start from the
same fiber are exponentially
approaching each other in forward time; similarly conclusion also holds for stochastic system \eqref{Equation-s}-\eqref{Equation-f}. \par
\noindent(iii) Stable foliation has geometric properties:   fibers of the stable foliation parallel to each other and the stable manifold is a special stable fiber.
\end{Theorem}
\par
\begin{Remark}\label{R4.6}
As $\mathbb{R}^{n+m}$ is a finite dimensional space, we can simply reserve the time to get the stable foliation by using the results of unstable foliation (Theorem \ref{Th4.1}, Theorem \ref{Th4.2}, Theorem \ref{Th4.3}). It is worth mentioning that different from the case of unstable foliation, the dynamical orbits which start from the
same stable fiber
approach each other in forward time in lower order $O(e^{p\gamma t})$ with $0<p<1$ rather than in $O(e^{\gamma t})$, but it does not affect the property of exponential approximation at all.
\end{Remark}

\renewcommand{\theequation}{\thesection.\arabic{equation}}
\setcounter{equation}{0}

\section{An example for unstable foliation}

\quad\; In this section, we present a simple   example for the theory developed in the previous section.
\par
Consider the following two dimensional   model with multiplicative L\'{e}vy noise in the framework of Marcus type SDEs
\begin{eqnarray}
\frac{dx}{dt}&=&x+x\diamond\dot{L_t^\alpha}, ~~~~~~~~~~\quad \hbox{in}\quad \mathbb{R}^1,\label{Equation-s-motivation}\\
\frac{dy}{dt}&=&-y+|x|+y\diamond\dot{L_t^\alpha}, ~ \quad \hbox{in}\quad \mathbb{R}^1,\label{Equation-f-motivation}
\end{eqnarray}
where $x$ (resp. $y$) is the unstable (resp. stable) component, accordingly, $a=1$, $b=-1$, $K=1$, $f(x,y)=0$, $g(x,y)=|x|$.
 \par
\noindent From Section 3, we can convert this SDE system to the following random system
\begin{eqnarray}
\frac{d\hat{x}}{dt}&=&\hat{x}+\hat{x}z(\theta_t\omega),~~~~~~~~~~~\quad \hbox{in}\quad \mathbb{R}^1,\label{Random-s-motivation}\\
\frac{d\hat{y}}{dt}&=&-\hat{y}+|\hat{x}|+\hat{y}z(\theta_t\omega),~~\quad \hbox{in}\quad \mathbb{R}^1. \label{Random-f-motivation}
\end{eqnarray}
Taking the initial value $\hat{x}(0)=\hat{x}_0$ and $\hat{y}(0)=\hat{y}_0$, we find the solution
\begin{eqnarray}
\hat{x}(t)&=&\hat{x}_0 e^{t+\int_0^tz(\theta_\tau\omega)d\tau},~t\in\mathbb{R},\\
\hat{y}(t)&=&\hat{y}_0 e^{-t+\int_0^tz(\theta_\tau\omega)d\tau}+\frac{1}{2}|\hat{x}_0|(e^{t+\int_0^tz(\theta_\tau\omega)d\tau}-e^{-t+\int_0^tz(\theta_\tau\omega)d\tau}),\quad t\in \mathbb{R},\notag\\
\end{eqnarray}
where
$$
\begin{array}{l}
\hat{x}(t)=\hat{x}(t,\omega, (\hat{x}_0, \hat{y}_0))
=\hat{x}(t,z(\theta_t\omega), (\hat{x}_0, \hat{y}_0))=x(t)e^{-z(\theta_t\omega)},\\
\hat{y}(t)=\hat{y}(t,\omega, (\hat{x}_0, \hat{y}_0))
=\hat{y}(t,z(\theta_t\omega), (\hat{x}_0, \hat{y}_0))=y(t)e^{-z(\theta_t\omega)},
\end{array}
$$
and $z(\theta_t\omega)=\int_{-\infty}^te^{(t-s)}dL_t^\alpha$ with the properties described in Lemma \ref{L3.1}.\par
\noindent On the one hand, it follows from Theorem \ref{Th4.1} that an unstable fiber of this system is described by \begin{equation}\label{foliation-motivation}
\mathcal{W}((\hat{x}_0, \hat{y}_0), \omega)=\{(\xi,
l(\xi, (\hat{x}_0, \hat{y}_0), \omega)) |\; \xi\in
\mathbb{R}^1\},
\end{equation}
where
\begin{eqnarray}\label{foliation-l-motivation}
l(\xi, (\hat{x}_0, \hat{y}_0),\omega)&=&\hat{y}_0+\int_{-\infty}^0e^{s+\int_s^0z(\theta_\tau\omega)d\tau}\cdot(|\xi|-|\hat{x}_0| ) e^{s+\int_0^sz(\theta_\tau\omega)d\tau}ds\notag\\
&=&\hat{y}_0+\frac{1}{2}(|\xi|-|\hat{x}_0|),\quad \xi\in
\mathbb{R}^1.
\end{eqnarray}
\noindent Moreover, by using the integral expression of the solution, for any two points $(\hat{x}_0, \hat{y}_0)$ and $(\hat{x}_0^\ast,\hat{y}_0^\ast)$ in $\mathbb{R}^1\times\mathbb{R}^1$, we calculate the difference between two orbits
$$
\begin{array}{lll}
J&:=|(\hat{x}(t,\omega, (\hat{x}_0, \hat{y}_0)),
\hat{y}(t,\omega, (\hat{x}_0, \hat{y}_0)))-(\hat{x}(t,\omega, (\hat{x}_0^\ast,\hat{y}_0^\ast)),\hat{y}(t,\omega, (\hat{x}_0^\ast,\hat{y}_0^\ast)))|\\
&=|\hat{x}(t,\omega, (\hat{x}_0, \hat{y}_0))-\hat{x}(t,\omega, (\hat{x}_0^\ast,\hat{y}_0^\ast))|+|\hat{y}(t,\omega, (\hat{x}_0, \hat{y}_0))-\hat{y}(t,\omega, (\hat{x}_0^\ast,\hat{y}_0^\ast))|\\
&\leq |\hat{x}_0-\hat{x}_0^\ast|\cdot e^{t+\int_0^tz(\theta_\tau\omega)d\tau}
+\frac{1}{2}|(\hat{x}_0-\hat{x}_0^\ast)|\cdot e^{t+\int_0^tz(\theta_\tau\omega)d\tau}\\
&~~~+|(\hat{y}_0-\hat{y}_0^\ast) -\frac{1}{2}(|\hat{x}_0|-|\hat{x}_0^\ast|)|\cdot e^{-t+\int_0^tz(\theta_\tau\omega)d\tau}.
\end{array}
$$
Recall that
$$\lim_{t \to \pm\infty}\frac{1}{t}\int_{0}^tz(\theta_\tau\omega)d\tau=0, \quad\omega \in \Omega$$i.e., $\int_{0}^tz(\theta_\tau\omega)d\tau$ has a sublinear growth rate, thus, the linear part of the exponent part in the exponential function plays a leading role.\par
\noindent Hence, if the coefficient
\begin{equation}\label{condition}
(\hat{y}_0-\hat{y}_0^\ast) -\frac{1}{2}(|\hat{x}_0|-|\hat{x}_0^\ast|)=0,
\end{equation}
then the difference of two orbits is
$J=O(e^{t})$, as $t\to -\infty$.
\par
\noindent We can obtain the following function
\begin{equation}\label{foliation-L-motivation}
L(\zeta, (\hat{x}_0, \hat{y}_0),
\omega)=\hat{y}_0+\frac{1}{2}(|\zeta|-|\hat{x}_0|),\quad \xi\in
\mathbb{R}^1,
\end{equation}
\par
\noindent which is in accordance with the function \eqref{foliation-l-motivation}, i.e., $l(\xi, (\hat{x}_0, \hat{y}_0),\omega)$. This immediately implies that the different dynamical orbits starting from the same fiber will be exponentially approaching each
other as $t\to -\infty$.
As seen in \eqref{foliation-l-motivation}, the unstable foliation of
\eqref{Random-s-motivation}-\eqref{Random-f-motivation} is a family
of the parallel curves (i.e., fibers) in the state space.
\par

In addition, from \eqref{manifold} and \eqref{manifold-h}, we see
that the \emph{unstable manifold} of
\eqref{Random-s-motivation}-\eqref{Random-f-motivation} is
\begin{equation}\label{manifold-motivation}
\mathcal{M}^{u}(\omega)=\{(\xi,
 h(\xi,\omega)) |\; \xi\in \mathbb{R}^1\},
\end{equation}
where
\begin{equation}\label{manifold-h-motivation}
h(\xi,\omega)=\frac{1}{2}|\xi|,\quad
\xi\in \mathbb{R}^1.
\end{equation}
By comparing with \eqref{foliation-l-motivation}, it is clear that the unstable manifold is a fiber of the unstable
foliation.
\par

\bigskip

\subsection*{Acknowledgements}
This work was partly supported by the NSFC grants 11531006 and 11771449.

The authors are grateful to Xianming Liu, Hongbo Fu and Ziying He for helpful discussions on stochastic equations driven by multiplicative L\'{e}vy noise.

\end{document}